\newbox\QEDbox
\newbox\nichts\setbox\nichts=\vbox{}\wd\nichts=2mm\ht\nichts=2mm
\setbox\QEDbox=\hbox{\vrule\vbox{\hrule\copy\nichts\hrule}\vrule}
\def\qed{\leavevmode\unskip\hfil\null\nobreak\hfill\copy\QEDbox\medbreak}

\documentclass{article}
\usepackage{amsmath}
\usepackage{amscd}
\usepackage{amsthm}
\usepackage{amsfonts}
\usepackage{amssymb}
\usepackage{verbatim}
\usepackage{epsfig}
\usepackage{mathtools}
\usepackage{color}
\usepackage{tikz}
\usepackage[all]{xy}
\usepackage{bm}
\renewcommand{\Bbb}{\mathbb}

\newcommand{\bu}{{\bf u}}
\newcommand{\bv}{{\bf v}}
\newcommand{\bw}{{\bf w}}

\newcommand{\bcr}{{\bf minA}}
\newcommand{\bs}{{\bf s}}

\newcommand{\bfm}{{\bf m}}

\newcommand{\BN}{\Bbb{N}}

\newcommand{\forb}[2]{\langle #1 \rangle_{#2}}
\newcommand{\xorb}[2][x]{\langle #1 \rangle_{#2}}
\newcommand{\orb}[1]{\forb{#1}{{}}}

\newcommand{\mpl}[1]{\ \mbox{mod}{}_+ {#1}}
\newcommand{\gws}[1]{[#1]}
\newcommand{\gw}[2]{[#1 \cdots #2]}
\newcommand{\mm}[2]{{\bfm}({#1,#2})}

\makeatletter
\newcommand{\superimpose}[2]{%
{\ooalign{$#1\@firstoftwo#2$\cr\hfil$#1\@secondoftwo#2$\hfil\cr}}}
\makeatother
\newcommand{\act}[1]{\ \mathpalette\superimpose{{\bigcup}{\mbox{\scriptsize $#1$}}}\ }

\theoremstyle{plain}
\newtheorem{Lemma}{Lemma}[section]
\newtheorem{Satz}{Theorem}[section]
\newtheorem{Cor}{Corollary}[section]

\theoremstyle{definition}
\newtheorem{Def}{Definition}[section]
\newtheorem{Rem}{Remarks}[section]
\newtheorem{Ex}{Example}[section]
\parskip \baselineskip
\begin{document}
\title{On infinite words avoiding a finite set of squares}
\date{October 20, 2013}
\author{Yasmine B. Sanderson}
\maketitle
\begin{abstract}
Building an infinite square-free word by appending one letter at a time while simultaneously avoiding the creation of squares is most likely to fail. When the alphabet has two letters this approach is impossible. When the alphabet has three or more letters, one will most probably create a word in which the addition of any letter invariably creates a square. When one restricts the set of undesired squares to a finite one, this can be possible. We study the constraints on the alphabet and the set of squares which permit this approach to work.
\end{abstract}
\large
\section{Introduction}
Suppose that one wanted to produce a (possibly infinite) sequence $w$ of symbols, or {\it word}, which would avoid a certain set $B$ of words. In other words, no element, or word, of $B$ should appear as a subsequence of $w$. The naive approach would be to build $w$ one symbol at a time, taking care at each step to avoid creating some element of $B$. It is well known that this approach does not always work. For example, if the set of symbols $A =\{a,b\}$ has only two elements ({\it letters}) and $B$ is the set of {\it squares}, that is, words of the form $vv$ where $v$ is a non-empty word on $A$, then this approach is thwarted after only a few steps: ajoining $a$ or $b$ to the word $bab$ produces either the square $bb$ or the square $abab$. We can think of $aba$ then as a {\it dead-end}; tagging on any letter from $A$ produces a forbidden subword.

When $A$ has three or more letters, the naive approach can work. However, the probability of obtaining such a word using this method is highly unlikely~\cite{berstel2005growth},~\cite{brandenburg1988uniformly}. In other words, one will most probably run into a dead-end at some point. Up to now, the main tool in producing squarefree words has been a squarefree morphism. Originally used by Thue~\cite{thue1912gegenseitige} in his computation of the first infinite squarefree word, variations have been used to construct squarefree words or words which avoid most squares, notably by Bean, Ehrenfeucht and McNulty~\cite{bean1979avoidable}, Berstel~\cite{Berstel1980235}, Carpi~\cite{Carpi1983231}, Crochemore~\cite{Crochemore1983235}, Dejean~\cite{Dejean197290}, Shallitt~\cite{shallit2004simultaneous}. Such a morphism maps squarefree words onto longer squarefree words. The infinite squarefree word is then obtained as the infinitely repeated iteration of this mapping. Some excellent surveys on the vast body of research in this area are given by Berstel~\cite{berstel2005growth},~\cite{berstel1984some},~\cite{Berstel2007996}, Currie~\cite{currie1993open},~\cite{Currie20057}, and Lothaire~\cite{lothaire1997combinatorics}.

When using these squarefree morphisms, creating the next part of the word requires repeating a recursive process. In addition, the process is deterministic; once a morphism and the initial word are chosen, there is no option to append a different letter. The need for such an approach is due to the number of (unwanted) squares being infinite. One can then ask whether the naive approach, or what we'll refer to from now on as the {\it sequential method}, works if the set of squares to be avoided is finite. Fix $\bs$, a finite sequence of natural numbers. Let $B(\bs)$ denote all words $vv$ in $A$ with the length $|v|\in \bs$. When $\bs = (1,2)$, then the above example shows that the answer is ``no'' for $A = \{a,b\}$. However, when $A = \{a,b,c\}$, then it works; there are no dead-ends because there are enough letters to prevent any potential squares. However, as soon as, say, $\bs = (1,3,5)$, then the sequential method no longer always works: affixing $a$, $b$, or $c$ to $cbacacbac$ produces squares with $v = cbaca$, $acb$ and $c$ respectively. In other words, with the sequential method we could come to a dead-end, namely $cabcacabc$.

In this paper we determine necessary and/or sufficient conditions on $|A|$ and $\bs$ for the existence of dead-ends for the sequential method. When there are no dead-ends then every word which avoids $B(\bs)$ is a subword in some (two-way) infinite word on $A$ which avoids $B(\bs)$. On the other hand, when dead-ends exist, we can not use the sequential method to build our infinite $B(\bs)$-avoiding word. 
We prove:
\begin{Satz}\label{thm-main} Let $A$ be an alphabet on $l >1$ letters. Let $\bs = (i_r, \ldots, i_2, i_1) \in \BN_{>0}^r$ be a strictly increasing sequence.

{\bf 1)} If $r<l$ then the sequential method has no dead-ends.

{\bf 2)} If $r=l$, $i_r \geq 2$, then the sequential method has dead-ends if and only if  either $\bs = ( i_r,  2i_r, \ldots, 2^{r-1}i_r)$ or there is a $0 < t < r-1$ such that $i_t \leq i_{t+1} + \cdots + i_{r}$. When $i_r=1$, these conditions are sufficient.

{\bf 3)}  If $ r> l$ then the sequential method has no dead-ends if every subsequence $\bs' = (i_{j_l},  \ldots, i_{j_2}, i_{j_1})$ $ \subset \bs$ of length $l$ is either of the form $(i_{j_l}, 2i_{j_l}, \ldots, 2^{l-1}i_{j_l})$ or satisfies $i_{j_t} \leq i_{j_{t+1}} + \cdots + i_{j_l}$ for some $0 < t < l-1$.
\end{Satz}
We conjecture that for $i_r\geq2$, the conditions {\bf 3)} are actually necessary. In Remark~\ref{rem-ir1} we discuss how the conditions {\bf 2-3)} are sufficient but not necessary for $i_r=1$.

The motivation for this problem comes from the development of products such as ``smart" random playback in music players~\cite{Leong:2006:RRD:1142405.1142428}, in which the order in which songs are played is not truly (i.e stochastically) random, but ``tweaked'' in order to give the impression of randomness. It has been shown~\cite{FK1} that humans rate sequences of symbols as being more likely to have been produced by a random process if the sequence avoids certain regularities, in particular squares~\cite{GT3}. Since human short-term memory is rather limited, it suffices to merely avoid a small set of them. In the case of the music playback, this paper addresses the question: if the playback uses the sequential method for song selection, how large does the playlist $|A|$ need to be in order to avoid unwished-for repeats $B(\bs)$.

We thank the Department of Cognitive Sciences, Central European University (Budapest), whose invitation provided valuable research time for this paper.

\section{Setup}

In this paper, we will use two different characterizations of words: one as set partitions and the other as strings of symbols. The difference in the two characterizations is mainly psychological, so we will frequently switch from one characterization to the one, choosing the one which facilitates understanding.

 For $n \in \BN$ let $\Pi_n$ be the set of partions of the set $[n]:=\{1,2,\ldots,n \}$. An element $\pi \in \Pi_n$ is then a set of {\it blocks } $\pi = \{ B_1, B_2, \ldots, B_i\}$ where each $B_t \not= \emptyset$, $B_s \cap B_t = \emptyset$ for each pair $1 \leq s \not= t \leq i$ and $\cup_{t=1}^iB_t = [n]$. The set $\Pi_n$ is equipped with a poset construction:
$$\pi \preceq \pi':=\{B_1', \ldots, B_k'\} \ \Leftrightarrow \forall s, \exists t \mbox{ such that } B_s \subseteq B_t' .$$
With this partial ordering we have a minimal element $\{ \{1\}, \{2\}, \ldots, \{n\}\}$ and a maximal element $\{ [n] \}$.
It is also graded $\Pi_n = \cup_{i= 1}^n\Pi_{n,i}$ with respect to the number of blocks. We denote by $k \sim j$ ($ k,j \in [n]$) when both $k$ and $j$ belong to the same block and $k \not\sim j$ when they belong to different blocks.

We recall that, if $A = \{ a_1, a_2, \ldots, a_l\}$ is an alphabet, then $$A^* = \{ w_1w_2\cdots w_m\ |\ m \in \BN_{\geq 0}, \ w_1, \ldots, w_m \in A\}$$ is the set of all {\it words} on $A$. For a word $w=w_1w_2\cdots w_m$,  $m$ is the {\it length} $|w|$ of $w$. 
We can represent each set partition $\pi = \{ B_1, B_2, \ldots, B_l\} \in \Pi_n$ by a certain string of symbols:
\begin{Def} A {\it generic word} ${\bf w}(\pi)$ of a set partition $\pi$ of $[n]$ is the sequence $w_1w_2\cdots w_n$ on $|\pi|$ letters which satisfies $w_k = w_j \Leftrightarrow k \sim j$ for all $1 \leq k,j \leq n$.\end{Def} Each generic word is then simply an equivalence class of words in sequence form. For example, the set partition $\{ \{1,3,4\}, \{2,6\}, \{5\} \}$ represents the word $w_2w_3w_1w_1w_2w_1$ where $w_1,w_2,w_3$ are any three different letters in an alphabet.
Since in this paper we are only interested in the property of whether $w_k = w_j$ and not which letter is associated to $w_k$ (or $w_j$), we will work interchangeably with set partitions and their associated generic words.

Let $\Sigma_n$ be the permutation group on $n$ elements. A permutation $\sigma \in \Sigma_n$ acts on a set partition $\pi = \{ B_1, B_2, \ldots, B_i\} \in \Pi_n$ in the obvious way: for each $B_t \in \pi$, $B_t = \{ i_1, \ldots, i_j \}$, we have $\sigma(B_t) := \{ \sigma(i_1), \ldots, \sigma(i_j) \}$. Then we define$$\sigma(\pi) := \{ \sigma(B_1), \ldots, \sigma(B_i) \}.$$ The permutation $\sigma$ also acts on the associated generic word in the obvious way:
$$\sigma({\bf w}(\pi)) = w_{\sigma(1)}w_{\sigma(2)} \cdots w_{\sigma(n)}.$$
A set permutation $\pi$ (resp. its generic word ${\bf w}(\pi)$) is $\sigma$-invariant if $\sigma(\pi) = \pi$ (resp. $\sigma({\bf w}(\pi)) = {\bf w}(\pi)$).

 We say that $w$ is {\it squarefree } if $w$ cannot be written as $w = xvvy$ with $|v| \geq 1$. For $i \in \BN_{>0}$, we define $w$ as being $i${\it -squarefree} if  $w$ cannot be written as $w = xvvy$ with $|v|=i$. For $\bs = (i_r, i_{r-1}, \ldots, i_1) \in \BN^r$, we define  $w$ as being $\bs$-{\it squarefree} if $w$ is $i_t$-squarefree for every $i_t \in \bs$. We say that $w$ has an {\it $i$-square starting at $j$} if we can write $w = xvvy$ where $|v|=i$ and $|x|=j-1$.

In the context of set partitions, $\pi$ has an $i$-square starting at $j$ if $j+t \sim j+t+i$ for $0\leq t \leq i-1$. It is $i${\it -squarefree} when there is no such $i$-square and $\bs$-{\it squarefree} when it is $i_t$-squarefree for all $i_t \in \bs$. Clearly, a word $w$ is $i$-squarefree if and only if $\pi_w$ is $i$-squarefree.
 
For a strictly increasing sequence $\bs \in\BN^r$, $\bs = (i_r < i_{r-1}< \cdots < i_1)$, let $B(\bs)$ be the set of all corresponding ``bad words''~\cite{noonan1999goulden}, in this case all square words $w=vv$ where $|v| = i_t$ for some $i_t \in \bs$. Then the set of all $\bs$-squarefree words is $A^* \backslash (A^*B(\bs)A^*)$. Let $N:=2i_1$ and let $G(\bs)$ to be the directed graph with vertex set $(A^* \backslash (A^*B(\bs)A^*))_N$.  Arcs, i.e. directed edges, are defined as follows: let $w = w_1w_{2}\cdots w_{N-1}w_N$ and $w'$ two vertices in $G(\bs)$. Then 
$$w \rightarrow w' \ \ \Leftrightarrow \ \ w' = aw_{1}\cdots w_{N-2}w_{N-1} \ \ \mbox{ for some } a \in A.$$

Every two-sided infinite $\bs$-squarefree words can be represented as an infinite path in $G(\bs)$ (and vice-versa). Each vertex represents the subword consisting of the $N$ last added letters. Each arc represents the last letter adjoined to the word. In order to guarantee that every random walk in this graph is infinite, there can be no dead-ends and no dead-starts as defined:

\begin{Def} A word $w \in G(\bs)$ is a {\it dead-end}, resp. {\it dead-start}, if the outdegree($w$)$=0$, resp. indegree($w$) $ = 0$.
\end{Def}

\section{ Characterization of dead-ends}

Let $\bs = (i_r, \ldots, i_1) \in \BN_{\geq 1}^r$ be strictly increasing and $N:=2i_1$. Then $w = w_1w_{2} \cdots w_{N-1}w_N \in G(\bs) = (A^*\backslash (A^*B(\bs)A^*))_N$ is a dead-end if and only if
$$ a \in A \ \ \ \Rightarrow \ \ \ w_2 \cdots w_{N-1}w_Na \in A^*B(\bs).$$
Analogously, $w$ is a dead-start if for every $ a \in A$, we have $aw_1w_{2}\cdots w_{N-1} \in B(\bs)A^*$. Since the reverse ordering $v_i\cdots v_1 v_i \cdots v_1$ of a square $v_1\cdots v_i v_1 \cdots v_i$ is also a square, we have that $w_Nw_{N-1} \cdots w_1$ is a dead-end if and only if $w_1w_{2} \cdots w_N$ is a dead-start. For ease of notation, we will actually determine the conditions under which dead-starts occur.

\medskip
From the fact that $ww = av = bv \ \ \Rightarrow \ \ a=b$ (that is, tagging on two different letters can not create squares of the same length), we obtain a first criterion:

\begin{Lemma}\label{lem-first}  If $|A| =l > r = |\bs|$ then $G(\bs)$ has no dead-ends. 
\end{Lemma}
In other words, when $|A| > |\bs|$, every word in $A^* \backslash (A^*B(\bs)A^*)$ is a subword of some infinite $B(\bs)$-avoiding word, which we can represent as some random walk on $G(\bs)$.

\bigskip
 The conditions $w_1\cdots w_N \in A^* \backslash (A^*B(\bs)A^*)$ and $ Aw_{1}w_2 \cdots w_{N-1}  \subseteq B(\bs)A^*$ puts restrictions on the letters $w_1, \ldots, w_{N-1}$: 

\begin{description}
\item[R1] $w_t = w_{t+i_j}$ for $0< t < i_j$, $1\leq j \leq r$, so that for some $a$, $aw_1\ldots w_{i_j-1} = w_{i_j}\ldots w_{2i_j}$.
\item[R2] $w_{i_j} = a_j$ for some unique $a_j \in A$. This forces $w_{i_j} \not= w_{i_k}$ for any two $i_j$, $i_k \in \bs$. 
\item[R3] That $w$ is $\bs$-squarefree to begin with forces $w_{i_j} \not= w_{2i_j}$.
\end{description}
We now look at these conditions in the context of set partitions. For each $i \in \BN_{n \geq 2}$ we define $\tau_i \in \Sigma_N$:
$$\tau_{i}(t)  := \begin{cases} t+i & \text{ for } 0 < t < i \\ t-i & \text{ for } i<t<2i \\ t & \text{ for } t=i \text{ or } 2i \leq t \leq N \\ \end{cases}.$$
 Condition R1 is equivalent to: $ \tau_{i_j}(\pi)=\pi$ for every $i_j \in \bs$. Conditions R2 and R3 are: $i_j \not\sim i_k$ for $j\not= k$ and $i_j \not\sim  2i_j$ for every $j$. The invariance condition R1 groups the elements of $[N]$ into a set of (disjoint) orbits which we denote by ${\bf o}(\bs)$. It is the minimal set partition $\pi$ which is invariant $ \tau_{i_j}(\pi) = \pi$ under the action of all $\tau_{i_j}$, $i_j \in \bs$. We define:
\begin{Def} The {\it primary difference conditions} on ${\bf o}(\bs)$ are 1) $i_j \not\sim i_k$ for $i_j,i_k \in \bs$, $i_j\not= i_k$, and 2) $i_j \not\sim 2i_j$ for every $i_j \in \bs$. Equivalently, the primary difference conditions on ${\bf w}(\bs)$ are 1) $\gws{i_j} \not= \gws{i_k}$ for $i_j,i_k \in \bs$, $i_j\not= i_k$, and 2) $\gws{i_j} \not= \gws{2i_j}$ for every $i_j \in \bs$.
\end{Def}

\begin{Ex}\label{ex-1a}  Let $\bs = (i) \in \BN_{\geq 1}$.  When $i \geq 2$, the action of $\tau_i$ gives $\xorb{i} = \{x,x+i\}$ for $0<x<i$ and $\xorb{i}=\{x\}$ for $x=i, 2i$. For $i=1$, ${\bf o}({\bs}) = \{ 1,2\}$. So
$${\bf o}(\bs) = \{ \forb{1}{i}, \ldots, \forb{i-1}{i}, \forb{i}{i}, \forb{2i}{i} \} \mbox{ with }\forb{x}{i} = \{ x, x+i \} \mbox{ for } 0 < x < i.$$
\end{Ex}

In order to determine ${\bf o}(\bs)$ for general $\bs$, we consider the (recursive) action of $\tau_{i_j}$ on ${\bf o}(i_{j-1},i_{j-2},\ldots, i_1)$. By $\forb{x}{i_{j-1}}$, we denote the orbit in ${\bf o}(i_{j-1},i_{j-2},\ldots, i_1)$ with minimal element $x$. The subscript $i_j$ represents the last action on it. (We use no subscript on the orbit $\orb{2i_1}$ since it is fixed by every $\tau_{i_j}$ for $i_j \in \bs$.) We use the notation  $$\forb{x}{i_{j-1}} \act{i_j} \forb{y}{i_{j-1}}$$
to denote that the action by $\tau_{i_j}$ conjoins the two orbits $\forb{x}{i_{j-1}}$ and $\forb{y}{i_{j-1}}$. 

\begin{Ex}\label{ex-2a} Let $\bs = (i_2, i_1) \in \BN^2_{\geq 2}$, $i_2 < i_1$. We consider the action of $\tau_{i_2}$ on ${\bf o}(i_1)$: for $0<x<\min(i_1-i_2,i_2)$, $\tau_{i_2}(x) = x+i_2 < \min(i_1,2i_2)$. In addition, $x+i_2$ is the minimal element in its orbit in ${\bf o}(i_1)$. So $\xorb{i_1} \act{i_2} \forb{x+i_2}{i_1}$. In the case that $i_1 - i_2 < i_2$, we have $\forb{i_1-i_2}{i_1} \act{i_2} \forb{i_1}{i_1}$. If the condition $i_1-i_2 < x < i_2$ is not empty, we have $\tau_{i_2}(x) = x+i_2$, but $x+i_2 > i_1$ so $x+i_2 \in \forb{x+i_2-i_1}{i_1}$ for all such $x$. This gives: $\xorb{i_1} \act{i_2} \forb{x+i_2-i_1}{i_1}$. With a change of variable in this last case, we can summarize the action of $\tau_{i_2}$ on ${\bf o }{(i_1)}$:
\begin{align}
\xorb{i_1} \act{i_2} \forb{x+i_2}{i_1} & \mbox{ for } 0<x<\min(i_1-i_2,i_2), \\
\forb{i_1-i_2}{i_1} \act{i_2} \forb{i_1}{i_1} & \mbox{ if } i_1 - i_2 < i_2, \\
\xorb{i_1} \act{i_2} \forb{x+i_1-i_2}{i_1} & \mbox{ for } 0 < x < \min(2i_2-i_1, i_2)
\end{align}
In other words, for $0 < x < \min(i_1-i_2, i_2)$ and for $ x = i_1-i_2$ (if $i_1-i_2< i_2 $), $$ \forb{x}{i_2} = \forb{x+i_2}{i_1} \bigcup^{m^*}_{m = 0} \forb{x+m(i_1-i_2)}{i_1} \mbox{ where } m^*:=\max \{ m \ | x + m(i_1- i_2) < i_2\}.$$
For $x = i_2, 2i_1$ and, when existant, $2i_2 \leq x \leq i_1$, the orbits $\xorb{i_1} \in {\bf o}(i_1)$ are fixed by $\tau_{i_2}$ because all elements in these orbits are either $= i_2$ or $\geq 2i_2$. So  $\xorb{i_2} =  \xorb{i_1}$.
\end{Ex}
We now determine ${\bf w}(\bs):= {\bf w}({\bf o}(\bs))$ for the above cases.
From here on, we use $\gws{t}$, resp. $\gw{s}{t}$ to denote the letter $w_t$, resp. subword $w_sw_{s+1}\cdots w_t$, in ${\bf w}(\bs)$. 
\begin{Ex}\label{ex-1b} From Example~\ref{ex-1a}, we see that ${\bf w}(i) = {\gw{1}{i-1}}\gws{i}{\gw{1}{ i-1}}\gws{2i}$ for $i \geq 2$. For $i=1$, ${\bf w}(1) = \gws{1}\gws{2}$. Clearly ${\bf w}(i)$ is $i$-squarefree if and only if $\gws{i}\not= \gws{2i}$. 
\end{Ex}
\begin{Ex}
It follows from Example~\ref{ex-2a} that the generic word ${\bf w}(i_2,i_1)$ is obtained via the image of ${\bf w}(i_1)$ under the map:
 $$\gw{1}{i_1} \mapsto \begin{cases} \gw{1}{i_2-1}\gws{i_2}\gw{1}{i_2-1}\gw{2i_2}{i_1} & \text{ if } 2i_2\leq i_1 \\
\gw{1}{i_1-i_2}^p\gw{1}{q-1}\gws{i_2}\gw{1}{i_1-i_2}  & \text{ if } 2i_2 > i_1 \\
\end{cases}$$
where $p \geq 0$, $0<q \leq i_1-i_2$ are determined by a (slightly modified) Euclidean division $i_2 = p(i_1-i_2) + q$. This gives us ${\bf w}(i_2,i_1) =$
\begin{align} 
 & {\bf u}\gws{i_2}{\bf u}\gw{2i_2}{i_1}{\bf u}\gws{i_2}{\bf u}\gw{2i_2}{i_1-1}\gws{2 i_1} & \text{ for } 2i_2 < i_1\label{genwordr2b} \\
&  {\bf u} \gws{i_2}{\bf u}\gws{i_1}{\bf u}\gws{i_2}{\bf u}\gws{2 i_1}& \text{ for } 2i_2 = i_1\label{genwordr2a}  \\
 & ({\bf u}\gws{i_1-i_2})^p\gw{1}{q-1}\gws{i_2}({\bf u}\gws{i_1-i_2}^{p+1}\gw{1}{q-1}\gws{i_2}{\bf u}\gws{2i_1} &  \text{ for } 2i_2 > i_1\label{genwordr2c} 
\end{align}
where ${\bf u}:=\gw{1}{(\min(i_2,i_1-i_2) -1)}$.
\end{Ex}

For a partition $\pi \prec {\bf o}(\bs)$ with $|A|$ blocks, ${\bf w}(\pi)$ is a dead-end for $G(\bs)$ if it $\bs$-squarefree and satisfies the primary difference conditions. It follows then that the size of $A$ is also a factor as to whether ${\bf w}(\bs)$ gives rise to dead-ends.
\begin{Ex} The partition ${\bf o}(3,5) = \{ \{1,4,6,9 \}, \{2,5,7\}, \{3,8 \}, \{10 \} \}$ with generic word ${\bf w}(3,5)=\gws{1}\gws{2}\gws{3}\gws{1}\gws{2}\gws{1}\gws{2}\gws{3}\gws{10}$ is $(3,5)$-squarefree. In an alphabet on $3$ or more letters, no word $\pi \prec {\bf w}(3,5)$ with $|\pi| = |A|$ is a dead-end; there are more letters than there are square lengths. For $A = \{a,b\}$, we consider a smaller partition with two blocks: $\pi := \{ \{ 1,2,4,5,6,7,9\}, \{3, 8, 10\} \}$. Since $\pi \prec {\bf o}(3,5)$, it is also $\tau_{5}$- and $\tau_{3}$-invariant. Its generic word ${\bf w}(\pi) = \gws{1}\gws{1}\gws{2}\gws{1}\gws{1}\gws{1}\gws{1}\gws{2}\gws{1}\gws{2}$ is $(3,5)$-squarefree and therefore gives dead-ends for $G(\bs)$ when $A$ has two letters:  $\gws{1}\gws{1}\gws{1}\gws{2}\gws{1}\gws{1}\gws{1}\gws{1}\gws{2}\gws{1}\gws{2}$ is a $5$-square and $\gws{2}\gws{1}\gws{1}\gws{2}\gws{1}\gws{1}\gws{1}\gws{1}\gws{2}\gws{1}\gws{2}$ has a $3$-square starting at $1$. The actual dead-ends in $G(\bs)$ would then be the words $ aabaaaabab$ and $bbabbbbaba$ for $A=\{a,b\}$. Clearly for an alphabet on one letter, all one-block partitions $\pi \prec {\bf o}(\bs)$ would have squares and therefore not be dead-ends.
\end{Ex}
Keeping this in mind, we define:

\begin{Def} The {\it minimal alphabet size} $\bcr(\bs)$ for the sequence $\bs$ is defined as 

1) $\infty$ if ${\bf o}(\bs)$ is not a candidate for dead-ends. In other words, either ${\bf o}(\bs)$ not $\bs$-squarefree or it does not satisfy the primary difference conditions,

2) the number $k \in \BN$ such that there exists a partition $\pi \preceq {\bf o}(\bs)$ on $k$ blocks which is a candidate for a dead-end and such that every smaller partition $\pi'\prec {\bf o}(\bs)$ on $k-1$ blocks is not.
\end{Def}
The primary difference conditions force $\bcr(\bs) \geq |\bs|$. 
In order to determine $\bcr$, it is therefore important to determine which orbits must remain dissimilar in order to prevent squares. In the case $\bs = (i)$ from Examples~\ref{ex-1a},~\ref{ex-1b}, it is necessary and sufficient that $\gws{i} \not= \gws{2i}$ in order to prevent $i$-squares. We can picture this as a 2-colouring of the {\it difference graph}: 
\begin{tikzpicture}[dot/.style 2 args={circle,inner sep=1pt,fill=blue, name=#1, label=#2}]
  \node [dot={a}{$\gws{i}$}] (a) at (0,0) {} ;
\node [dot={b}{$\gws{2i}$}] (b) at (2,0) {} ;
  \draw  (a) -- (b);
\end{tikzpicture}
It follows that $\bcr(i) = 2$. 
\begin{Lemma}\label{lem-2b} Let $\bs = (i_2, i_1) \in \BN^2_{\geq 2}$, $i_2 < i_1$. If $2i_2=i_1$ then $\bcr(i_2,i_1)=3$. Otherwise $\bcr(i_2,i_1) = 2$ .
\end{Lemma}

{\bf Proof:} From Example~\ref{ex-1a} we know that ${\bf o}(\bs)$ is $i_1$-square free if and only if $\gws{i_1} \not= \gws{2i_1}$. We now determine conditions to avoid $i_2$-squares.

{\bf Case $2i_2 = i_1$: }  From (\ref{genwordr2a}) we see that ${\bf w}(\bs)$ is $\bs$-squarefree if and only if $\gws{i_2} \not= \gws{2i_2}=\gws{i_1}$, $\gws{i_2} \not= \gws{2i_1}$ and $\gws{i_1} \not= \gws{2i_1}$. The difference graph is the complete graph on three vertices: $i_1$, $i_2$, $2i_1$. It is $3$-colorable but not $2$-colorable. Therefore $\bcr(i_2,i_1)=3$.

{\bf Case $ 2i_2 < i_1$: } We refer to (\ref{genwordr2b}) for the form of ${\bf w}(i_2,i_1)$. 
Since $\gws{i_1}$ is a prefix of $\gw{1}{i_1}$ in ${\bf w}(i_1)$, it is also a prefix of $\gw{1}{i_2-1}\gws{i_2}\gw{1}{i_2-1}\gws{2i_2}$ in ${\bf w}(i_2,i_1)$. Avoiding $i_2$-squares requires then $\gws{i_2} \not= \gws{2i_2}$ and $\gws{i_2} \not= \gws{i_1}$. These, along with $\gws{i_1} \not= \gws{2i_1}$  are (necessary and) sufficient for $\bs$-squarefreeness only in the intervals $ w\gws{i_2}w\gws{2i_2}$ and $\gws{i_1}w\gws{i_2}w\gws{2i_2}$. In order to guarantee $i_2$-squarefreeness in all of ${\bf w}(\bs)$, one must set $\gws{x} \not= \gws{x+i_2}$ for enough $x$ such that every subword in ${\bf w}(\bs)$ of length $2i_2$ contains such a pair. There are possibly many ways of doing this. We give one: set
\begin{align}\gws{i_1-1} \not= \begin{cases} \gws{i_1-i_2-1} & \text{ if } i_1-i_2-1 \geq 2i_2 \\ 
 \gws{i_1-2i_2-1} & \text{ if } i_1-i_2-1 < 2i_2 
\end{cases}.\label{diff-1}\end{align}
This guarantees $i_2$-squarefreeness for every subword of length $2i_2$ and containing $\gws{i_1}$ in the concatenation: $$[ 2i_2 \cdots  \underbrace{s' \cdots i_1-1 ] \gws{ i_1} [ 1 \cdots s}_{\text{subword of length $2i_2$}} \cdots i_2-1 ].$$ If $i_1 - 2i_2 > i_2$ then additionally set \begin{align}\gws{2i_2} \not= \gws{3i_2} \not= \cdots = \gws{fi_2}, \mbox{ where } f:= \max \{ t \ | \ ti_2 < i_1 \}.\label{diff-2}\end{align}
These conditions do not require more than $2$ letters and the resulting difference graph is $2$-colorable. So $\bcr(i_2,i_1) = 2$.

{\bf Case $2i_2 > i_1$:} We refer to (\ref{genwordr2c}) for the form of ${\bf w}(i_2,i_1)$. In ${\bf w}(i_2,i_1)$ the letter $\gws{2i_2} = \gws{q}$ represents $2i_2$ and $i_1 + 2i_2$. The letter $\gws{i_1-i_2}$ represents $2i_1-i_2$, and $i_1$. A necessary condition for $i_2$-squarefreeness is then $\gws{i_2} \not= \gws{2i_2}$ and $\gws{i_2} \not= \gws{i_1-i_2}$. This is in fact sufficient: we have
$$\xymatrix@1@=1pt@M=0pt{ i_2 \ar@{<.>}@/^1pc/[rrrr]|{\ \not\sim \ }& <&  i_1 \ar@{<.>}@/_1pc/[rrrr]|{\ \not\sim \ }  &<&  2i_2 & <   i_1&+&i_2 & <  2i_1}  \ \ \mbox{ and } \xymatrix@1@=1pt@M=0pt{ 2i_1&-\ar@{<.>}@/_1pc/[rrrrr]|{\ \not\sim \ } &i_2 & <&  i_1+i_2 & < & 2i_1}.$$
Every subinterval in $0 < x< 2i_1$ of length $2i_2$ contains at least one of the dissimilar pairs: $\gws{i_2}\not= \gws{2i_2}$, or $\gws{i_1} \not= \gws{i_1+i_2} = \gws{i_2}$, or $\gws{2i_1-i_2} = \gws{i_1} \not=\gws{2i_1}$. Restriction R2 forces $\gws{i_1} \not= \gws{i_2}$. The resulting difference graph 
\begin{tikzpicture}[dot/.style={circle,inner sep=1pt,fill=blue}]  
\node (1) at (0,0) [dot, label=90:$\gws{2i_2}$] {};  
\node  (2) at (1,0) [dot,label=90:$\gws{i_2}$]  {};
 \node (3) at (2,0) [dot, label=90:$\gws{i_1}$]  {};
 \node (4) at (3,0) [dot, label=90:$\gws{2i_1}$]  {};
 \draw  (1) -- (2)
(2) -- (3)
(3) -- (4);
\end{tikzpicture} 
is clearly $2$-colorable, so $\bcr(i_2,i_1) = 2$.

\qed

\begin{Ex}\label{ex-3a} Let $\bs = (1, i_1) \in \BN^2_{\geq 1}$, $i_1 > 1$. For the generic word ${\bf w}(1, i_1)$ to give dead-ends, there can be no $1$-squares. This means that $\gws{x} \not= \gws{x+1}$ for any $0 < x < 2i_1$. The difference graph has the following form: the vertices are the orbits in ${\bf o}(i_1)$. An edge connects vertex $y_1$ with vertex $y_2$ when $x \in \forb{y_1}{i_1}$, $x+1 \in \forb{y_2}{i_1}$. Together with an edge representing $\gws{i_1} \not= \gws{2i_1}$ this gives us:
 $$
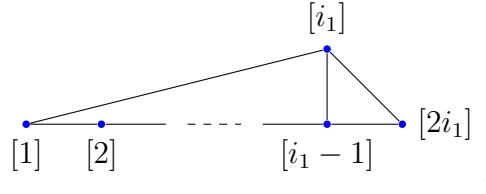
\begin{tikzpicture}[dot/.style={circle,inner sep=1pt,fill=blue}] 
 \node (1) at (0,0) [dot, label=-90:$\gws{1}$] {} ;
        \node (2) at  (1,0) [dot, label=-90:$\gws{2}$]  {} ;
      \node (i11) at (4,0) [dot, label=-90:$\gws{i_1-1}$] {} ;
     \node (i1) at  (4,1) [dot, label=90:$\gws{i_1}$] {} ;
   \node (2i1) at   (5,0) [dot, label=right:$\gws{2i_1}$] {} ;
    \node (3) at   (2,0) [] {};
    \node (4) at  (3,0) [] {};
  \draw  (1) -- (2)
(2) -- (3)
(4) -- (i11)
(i11) -- (i1)
 (2i1) -- (i1)
(i1) -- (1)
(i11) -- (2i1);
\begin{scope} [dashed]
\draw (3) -- (4) ;
\end{scope}
\end{tikzpicture},$$
a graph which is $3$- but not $2$-colorable, so $\bcr(1, i_1) = 3$. \end{Ex}
 When $\bs$ is of length two, then there are no dead-ends when $|A| >2$ or when any potential dead-end already has a square in it and is therefore not included as a vertex in $G(\bs)$. From Lemma~\ref{ex-2a} and Example~\ref{ex-3a} we obtain:
\begin{Cor} Let $\bs = (i_2, i_1) \in \BN^2_{\geq 1}$, $i_2 < i_1 $. Then $G(\bs)$ has no dead-ends iff either $|A| \geq 3$, or $2i_2=i_1$, or $i_2=1$.\end{Cor}

\section{General theorems}
The rest of this paper will be devoted to determining $\bcr(\bs)$ for $\bs$ of arbitrary length. We show that $\bcr(\bs) < \infty$  only if $\bs$ satisfies certain conditions (which we will call condition C). This condition is necessary and sufficient for the primary difference conditions and almost sufficient for $\bs$-squarefreeness.

\begin{Def} A sequence $\bs= (i_r,\ldots, i_2, i_1) \in \BN^r_{\geq 2}$ satisfies {\bf condition C} if the following inequalities hold: 
\begin{eqnarray*} i_1 &>& i_2 + \cdots + i_{r-1}+  i_r  \\
i_2 &>&  i_3 + \cdots  + i_{r} \\
\vdots &> & \vdots \\
i_{r-2} &> & i_{r-1}  + i_r \\
i_{r-1} &> &i_r
\end{eqnarray*}
\end{Def}

We will use the following notation. Let $\bs = (i_r,\ldots, i_2, i_1)$ be a sequence satisfying condition C. For $1 \leq v \leq r$ we define $\mm{v}{v} := i_v$, $\mm{v+1}{v}:= 2i_v$, and for $1 \leq u < v \leq r$,  $$\mm{u}{v}:=\min(i_u - \sum_{j=u+1}^{v}i_j,\  i_{u+1} -\sum_{j=u+2}^{v}i_j , \ldots,\  i_{v-1}-i_v,\  i_v).$$ 
 We note that \begin{align} \mm{u}{v}=\min(i_u - \sum_{j=u+1}^{v}i_j, \mm{u+1}{v}) \text{ and }\\ \mm{u}{v+1})= \min(\mm{u}{v}-i_{v+1},i_{v+1})\end{align} It follows that $\mm{u+1}{v}\geq \mm{u}{v} > \mm{u}{v+1}$ and \begin{align}\mm{v+1}{v+1}=i_{v+1} < \mm{1}{v} \leq \mm{2}{v} \leq \cdots \leq \mm{v-1}{v} \leq \mm{v}{v} = i_v \label{ineq1}.\end{align}
 Although $\mm{u}{v}$ depends on the choice of sequence $\bs$, we will omit $\bs$ from the notation because it will be always clear which sequence is being used.

We will need an explicit description of the orbits of ${\bf o}(\bs)$:
\begin{Satz}\label{thm-orbits} Let $\bs = (i_r,  \ldots, i_2, i_1) \in \BN_{\geq 2}^r$ satisfy condition C.
 
{\bf 1)} Then $\tau_{i_r}$ acts non-trivially on the orbits of ${\bf o}(i_r, \ldots, i_{1})$ in the following way:
\begin{align}
&\xorb{i_{r-1}} \act{i_r} \forb{x+i_r}{i_{r-1}} & 0 < x < \mm{1}{r}\label{case1-eq1} \\
&\forb{\mm{t}{r}}{i_{r-1}} \act{i_r} \forb{\mm{t}{r-1}}{i_{r-1}}  & \mbox{ when } \mm{t}{r} < \mm{r}{r}\label{case1-eq2} \\
& \xorb{i_{r-1}} \act{i_r}  \forb{x+ \mm{t}{r}}{i_{r-1}} & \mbox{ for } 0 < x < \mm{t+1}{r}-\mm{t}{r} \label{case1-eq3}
\end{align}
where $1 \leq t \leq r$.

{\bf 2)} For $0 < x < \mm{1}{r}$ and for $x = \mm{t}{r}$ for $t$ such that $\mm{t}{r} < \mm{r}{r}$,
\begin{align}\xorb{i_r} = \bigcup_{s= 1}^r \bigcup_{n_s \geq 0}^{n^*_s} 
\forb{x + n_s(\mm{s}{r})}{i_{r-1}} \label{case2}\end{align}
where  $n_s^*:=\max \{ n \in \BN \ | \ x + n(\mm{s}{r}) < \mm{s+1}{r}\}.$
\end{Satz}

\noindent
{\bf Proof:} The proof is by induction on $r$. We have already seen that both claims 1) and 2) hold for $r=2$. Now let $r > 2$. The operator $\tau_{i_r}$ acts non-trivially on orbits containing the elements $x$ with $0 < x < i_r$. Because $i_r < \mm{1}{r-1}$, every element $x$, $ 0 < x < i_r$,  lies in the orbit $ \xorb{i_{r-1}}$. Since $\tau_{i_r}(x) = x+i_r$ for $0 < x < i_r$, we have $\xorb{i_{r-1}} \act{i_r} \forb{x+i_r}{i_{r-1}}$ whenever $x + i_r$ is the representative (or smallest element) of its orbit. In other words, for $0 < x< \min(\mm{1}{r-1}-i_r,i_r) = \mm{1}{r}$. This gives us the action (\ref{case1-eq1}) and part of the equality (\ref{case2}). 

 When $x+i_r \geq \mm{1}{r-1}$, then there are two cases to consider. The first is when $x+i_r=\mm{t}{r-1} < 2i_r$ for some $t$. Then $x = \mm{t}{r} < i_r$. We thus obtain action (\ref{case1-eq2}): $$ \forb{\mm{t}{r-1}}{i_{r-1}} \act{i_r} \forb{\mm{t}{r}}{i_{r-1}}.$$

In the second case $\mm{t}{r-1} < x+i_r < \mm{t+1}{r-1} \leq 2i_r$ for some $t$. Since $i_r < \mm{t}{r-1}$ for all $t$, we have $2i_r< 2\mm{t}{r-1}$. By induction on the orbit description (\ref{case2}), it holds that $x+i_r = y + \mm{t}{r-1}$ for some $0 < y < \mm{1}{r-1}$. In other words, $\tau_{i_r}$ conjoins the orbits $\forb{x}{i_{r-1}}$ and $ \forb{y}{i_{r-1}} = \forb{x - \mm{t}{r}}{i_{r-1}}$. A change of variable gives us action (\ref{case1-eq3}) and the rest of the equality (\ref{case2}). 
 \qed

\medskip
\bigskip
\begin{Cor}~\label{cor2}
The generic word ${\bf w}(i_r,\ldots, i_1)$ is obtained from ${\bf w}(i_{r-1},\ldots, i_{1})$ by the map \begin{align}\gw{1}{\mm{1}{r-1}} \mapsto \begin{cases}\gw{1}{i_r}\gw{1}{i_r-1}\gw{2i_r}{\mm{1}{r-1}} & \text{ when } 2i_r \leq \mm{1}{r-1}\\
\gw{1}{i_r}\gw{1}{\mm{1}{r}} &  \text{ when } 2i_r > \mm{1}{r-1}\\
\end{cases}\label{action}\end{align}
where the generic word $\gw{1}{i_r}$ is obtained via successive applications of the recursive (in $t$) mappings:\begin{equation}\label{action2}\gw{1}{\mm{t+1}{r}} \mapsto \left(\gw{1}{\mm{t}{r}}\right)^{p_{t,r}}\gw{1}{q_{t,r}-1}\gws{\mm{t+1}{r}}, \ \ t=r-1,\ldots, 2,1.\end{equation}
where $\mm{t+1}{r} = p_{t,r}\mm{t}{r} + q_{t,r}$ with $ p_{t,r} \geq 0$ maximal and $0 <  q_{t,r} \leq \mm{t}{r}$.
\end{Cor}
{\bf Proof:} We have seen in Examples~\ref{ex-1a},~\ref{ex-2a} that the Corollary holds for $r\leq 2$. The action of $\tau_{i_r}$ results in $\forb{x}{i_{r-1}} \act{i_r} \forb{x+i_r}{i_{r-1}}$ for all $0 < x < i_r$. This gives the mapping~\ref{action}. For $\mm{1}{r} = i_r$, the mappings of type~\ref{action2} are then identity maps and therefore hold trivially. For $\mm{1}{r} = \mm{1}{r-1}-i_r < i_r$, we have $\mm{1}{r} \leq \mm{2}{r} \leq \cdots \leq \mm{r}{r}=i_r < \mm{1}{r}$. Fix $t$, $1\leq t \leq r-1$. For every $x$, $\mm{t}{r} < x \leq \mm{t+1}{r}$, we have $x = p_{t,r}(x)\mm{t}{r} + q_{t,r}(x)$ with $ p_{t,r}(x) \geq 0$ maximal and $0 <  q_{t,r}(x) \leq \mm{t}{r}$. By Theorem~\ref{thm-orbits}-2, $\gws{x} = \gws{ q_{t,r}(x)}$. This gives the mapping~\ref{action2}.
\qed

We now show that ${\bf o}(\bs)$ (and ${\bf w}(\bs)$) satisfies the primary difference conditions if $\bs$ satisfies Condition C.
\begin{Satz}\label{cor-2} Let $\bs \in \BN^r_{\geq 2}$ satisfy Condition C.

{\bf 1)} For every $1 \leq s \leq r$, there exists a $t$, with $s \leq t \leq r$, such that $i_s \in \forb{\mm{s}{t}}{i_r}$. It follows that the $i_1, \ldots, i_r$ are pairwise dissimilar.

{\bf 2)} For every $1 \leq t \leq r$, we have $i_t \not\sim 2i_t$.

{\bf 3)} There exists an $0 < x < 2i_1$ with $x \sim x+1$ iff $i_s-i_{s+1} -i_{s+2}-\cdots - i_t = 1$ for some $1 \leq s < t \leq r$.
\end{Satz}
{\bf Proof:}  From Theorem~\ref{thm-orbits}-1), we have that $i_s \in \forb{\mm{s}{s}}{i_t}=\forb{i_{s}}{i_t}$ for $t\leq s$. For $t > s$, we have $i_s \in \forb{\mm{s}{t}}{i_{t}}$ if and only if 
\begin{align*}i_s - i_{s+1} < i_{s+1} &\ \ \ \Rightarrow \forb{\mm{s}{s+1}}{i_{s+1}} \supseteq \forb{\mm{s}{s+1}}{i_{s}} \cup \forb{i_s}{i_{s}} \\ 
i_{s} - i_{s+1} - i_{s+2} < i_{s+2}&\ \ \ \Rightarrow \forb{\mm{s}{s+2}}{i_{s+2}} \supseteq \forb{\mm{s}{s+2}}{i_{s+1}} \cup  \forb{\mm{s}{s+1}}{i_{s+1}}  \\ \ldots & \ldots \\ i_s - i_{s+1} - \cdots - i_{t\ \ }  <i_t &\ \ \ \Rightarrow \forb{\mm{s}{t}}{i_{t}} \supseteq \forb{\mm{s}{t}}{i_{t-1}} \cup  \forb{\mm{s}{t-1}}{i_{t-1}} \end{align*} Equivalently, $i_s \in \forb{\mm{s}{t}}{i_{t}}$ if and only if $\mm{s}{t} < \mm{s+1}{t}$. By construction, we have $\forb{\mm{s}{t}}{i_t} \not= \forb{\mm{s'}{t'}}{i_t}$ for $\mm{s}{t} \not= \mm{s'}{t'}$ so the $i_s$ are pairwise dissimilar. This proves claim 1.

  We first show that the action of $\tau_{i_r}$ fixes $\forb{i_r}{i_{r-1}}$. By Theorem~\ref{thm-orbits}{\bf -2)}, we have
$$\forb{i_r}{i_{r-1}} = \cup_{t=1}^{r-1} \cup_{n_t \geq 0 } \forb{i_r+n_t(\mm{t}{r-1})}{i_{r-2}}$$ with $n_t$ bounded by the restriction $i_r + n_t(\mm{t}{r-1}) < \mm{t+1}{r-1}$. For $n_t>0$, the smallest element in $ \forb{i_r+n_t(\mm{t}{r-1})}{i_{r-2}}$ is $i_r+n_t(\mm{t}{r-1})$, which satisfies $i_r+n_t(\mm{t}{r-1}) \geq i_r + i_t - i_{t+1} - \cdots - i_{r-1} > 2i_r$. Likewise, the smallest element in $\forb{i_r + i_{r-1}}{i_{r-2}}$ is $i_r + i_{r-1}$ with $i_r + i_{r-1} >2i_r$. For $n_t=0$, the second smallest element is (for some $t$) $i_r + \mm{t}{r-2}>2i_r$. It follows that every element in $\forb{i_r}{i_{r-1}}$ is therefore either $= i_r$ or $> 2i_r$ and so is fixed by $\tau_{i_r}$. Since $2i_r \not\in \forb{i_r}{i_{r-1}}$, it follows that $i_r \not\sim 2i_r$ in ${\bf o}({\bs})$. The rest of the proof of claims 1-2 is by induction. The claim is shown for $r=1,2$ in Examples \ref{ex-1a} and \ref{ex-2a}. We assume that the claims hold for $r-1$. We need to show that they continue to hold under the action of $\tau_{i_r}$. In fact, we need only look at the cases where {\it both} $\forb{\mm{s}{t}}{i_{r-1}}$ and $\forb{x}{i_{r-1}} \ni 2i_s$ are {\it not} fixed by $\tau_{i_r}$ and make sure that they are not conjoined under the action of $\tau_{i_r}$. The orbit $\forb{\mm{s}{t}}{i_{r-1}}$ is not fixed by $\tau_{i_r}$ if and only if $\mm{s}{t}=\mm{s}{r}<\mm{s+1}{r} $. We then have $i_s \in \forb{\mm{s}{r}}{i_{r}}$. On the other hand, we have $2i_s \in \forb{2i_s - i_{s+1}}{i_{s+1}}$. Since $$2i_s - i_{s+1} < 2i_s - (i_s + i_{s+1} + \cdots + i_t) = i_s -i_{s+1} - \cdots -i_t \leq \mm{s}{t}$$ for all $s+1 \leq t \leq r$, then we have
$$2i_s \in  \forb{2i_s - i_{s+1}}{i_{s+1}} \subset \forb{x}{i_r} $$ for some $x < \mm{s}{r}$. By Theorem~\ref{thm-orbits}, $\forb{x}{i_r} \cap \forb{\mm{s}{r}}{i_r} = \emptyset$, so $i_s$ and $2i_s$ remain dissimilar. It follows that claim 2 hold for all $s$, $1 \leq s \leq r$.

Proof of {\bf 3)} {\bf ($\Leftarrow$)} Without loss of generality, we can assume $s=1$, $t=r$ and that $i_f-i_{f+1} - \cdots - i_g >1$ for every $1<f<g < r$.  The double inequality $$i_f < i_1 - i_2 - \cdots - i_{f-1} = i_f + \cdots i_r +1 \leq 2i_f$$ holds for every $1 < f \leq r$. Using the same argument as in the proof of 3) we see that the orbit $\forb{\mm{1}{r}}{i_r} = \forb{1}{i_r}$ has smallest element $\mm{1}{r}=1$ and contains $i_1$. From the action of $\tau_{i_1}$ it also contains $1+i_1$ so $i_1 \sim i_1+1$. {\bf ($\Rightarrow$)} For the other direction, if $y+1, y \in \forb{x}{i_r}$ then $y+1 = x + k_{t'}\mm{t'}{s'}$ and $y = x + k_{t}\mm{t}{s} < \mm{t+1}{s}$ for some $ 1< s' < s\leq r$, $ 1 \leq t \leq s+1$, $1 \leq t' \leq s'+1$. From inequality chain~\ref{ineq1}, we see that we have three cases to consider. The first is when $\mm{t}{s} = \mm{s}{s} = i_s$ and $\mm{t'}{s'} = \mm{1}{s-1} <2i_s$. Then $k_{t'} = k_t=1$ and $1 = y+1 -y = \mm{1}{s-1} - i_s = \mm{1}{s}$. For the next case, we have $y = x + k_{t}\mm{t}{s} < \mm{t+1}{s} < x + k_{t'}\mm{t'}{s'}$. Then $k_{t'}\mm{t'}{s'} - k_{t}\mm{t}{s} =1$ only if $k_t = 0$, $k_{t'}=1$ and $\mm{t'}{s'}=1$. In the third case, $\mm{t'}{s'} = \mm{t}{s}$. Then $k_{t'}=k_t+1$ and $\mm{t}{s}=1$.
\qed

\begin{Rem}\label{rmk-1} In general, the primary difference conditions are (almost) necessary for squarefreeness. Specifically, in order for $\bcr(\bs) < \infty$, it is necessary that $i_r \not\sim 2i_r$ and that $i_r \not\sim i_j$ for all $j< r$ with $2i_r < i_j$. In fact, since $\tau_{i_r}(t) = t+ i_r$ for all $0<t< i_r$, the additional condition $i_r \sim 2i_r$ gives $t \sim t+i_r$ for $1 \leq t\leq i_r$, that is, a $i_r$-square starting at $1$. Likewise, when $i_j + 2i_r < 2i_j$,  we have for $0<t<i_r < i_j$
$$ \tau_{i_j}\tau_{i_r}\tau_{i_j}(i_j + t) = \tau_{i_j}\tau_{i_r}(t) = \tau_{i_j}(t+i_r) = t+i_r+i_j.$$
The additional condition $i_r \sim i_j$($\sim i_j + i_r$) gives $t \sim t+i_r$ for $i_j \leq t \leq i_j+i_r-1$, that is, a $i_r$-square starting at $i_j$.
\end{Rem}
From Example~\ref{lem-2b} we know that the primary difference conditions are not necessarily sufficient for $\bs$-squarefreeness. Any additional difference conditions that we will impose in order to guarantee this will be defined as {\it secondary difference conditions}.
\begin{Satz}\label{thm-pita}  Let $\bs = (i_r, \ldots, i_2, i_1) \in \BN_{\geq 2}^r$ satisfy condition C. Then $$\bcr(\bs)= \begin{cases} r+1 & \text{ if } \bs = (i_r,  2i_r, \ldots,  2^{r-1}i_r) \\  r & \text{ otherwise } \end{cases}.$$
\end{Satz}

{\bf Proof of Theorem~\ref{thm-pita}:} Examples~\ref{ex-1a},~\ref{ex-2a} show that the claim holds for $r=1,2$.

 In the case that $\bs = ( 2^{r-1}i_r,  2^{r-2}i_r, \ldots, i_r)$, $${\bf w}(i_1, \ldots, i_j) = \prod_{k=j}^1\gw{1}{i_j-1}\gws{i_j}\gw{1}{i_j-1}\gws{\mm{k+1}{k}} = \prod_{k=j}^1\gw{1}{i_j-1}\gws{i_j}\gw{1}{i_j-1}\gws{2i_k}$$
Every instance of $\gw{1}{i_j-1}\gws{i_j}\gw{1}{i_j-1}$ has suffixes $2i_j=i_{j-1}$, $2i_{j-1}=i_{j-2}$, $\ldots$, $2i_2 = i_1$, $2i_1$ and prefixes $\gws{\ }$, $\gws{i_{j-1}}, \ldots, \gws{i_1}$. In order to prevent $i_j$-squares it is necessary and sufficient that $\gws{i_j} \not= \gws{i_k}$ for all $k \not= j$ and $\gws{i_j} \not= \gws{2i_1}$. This holds for all $1\leq j \leq r$. Therefore the difference graph is the complete graph on the vertices $\gws{i_1}, \ldots, \gws{i_r}, \gws{2i_1}$, which is $(r+1)$-colorable. 

 From here on we consider the case that at least one $i_j$ satisfies $2i_j \not= i_{j-1}$ and analyze to what extent the primary difference conditions ensure squarefreeness. The difference graph coming from these primary difference conditions contains the complete graph on the $r$ vertices $\gws{i_1}, \cdots, \gws{i_r}$. The condition $2i_j \not= i_{j-1}$ ensures that the difference graph is not complete on $r+1$ vertices and is therefore is $r$-colorable.
 We show that secondary difference conditions do not require more than $r$ letters. 

{\bf I - } We first consider the case where $\mm{1}{t} = \mm{1}{t-1}-i_t$ for all $1 <t\leq r$.
The action of $\tau_{i_t}$ replaces every subword $\gw{1}{\mm{1}{t-1}}$ in $\bw(i_{t-1},\ldots,i_{1})$ by $\gw{1}{i_t-1}\gws{i_t}\gw{1}{\mm{1}{t}}$, where $\gw{1}{i_t-1}$ is then mapped according to Corollary~\ref{cor2}. 
By induction, $\bw(i_{t-1},\ldots,i_2,i_1)$ is of the form
$$\bu_{1}\gws{\mm{s_1}{t_1}}\bu_{2}\gws{\mm{s_2}{t_2}}\cdots \bu_{j}\gws{2i_1}$$
where each $\bu_{j_k}$ is a (mixed) power of $\gw{1}{\mm{1}{t-1}}$:
$$\bu_{k} = \left(\gw{1}{\mm{1}{t-1}}\right)^{p_{k}}\gw{1}{q_{k}-1}.$$
  The image $\tau_{i_t}(\bu_k)\subset {\bf w}(i_t,\ldots, i_1)$ of each $\bu_k$ is $i_t$-squarefree: set $\bfm:=\mm{1}{t}$. Then in $\bu_{k}$ we have the following overlapping pairs:
$$\xymatrix@1@=1pt@M=0pt{
i_t\ar@{<.>}@/^1pc/[rrrr]|{\ \not\sim \ } &<& \bfm \ar@{<.>}@/_1pc/[rrrr]|{\ \not\sim \ } &<& 2i_t &<& i_t + \bfm \ar@{<.>}@/^1pc/[rrrr]|{\ \not\sim \ } &<& 2\bfm \ar@{<.>}@/_1pc/[rrrr]|{\ \not\sim \ }&<& 2i_t + \bfm &<& i_t + 2\bfm \ar@{<.>}@/^1pc/[rrrr]|{\ \not\sim \ }&<& 3\bfm \ar@{<.>}@/_1pc/[rrrr]|{\ \not\sim \ }&<& 2i_t + 2 \bfm &<&\cdots }$$

which in $\tau_{i_t}(\bu_{j_k})$ are assigned different letters: $\gws{i_t + n\mm{1}{t}} = \gws{i_t } \not= \gws{2i_t}= \gws{2i_t + n\mm{1}{t}}$, $\gws{n\mm{1}{t}}= \gws{i_1} \not= \gws{i_t} = \gws{i_t + n\mm{1}{r}}$ ($n \geq 0$). Therefore the $\bu_j$ themselves are $i_t$-squarefree. By Theorem~\ref{cor-2}, these non-equal letters remain different in ${\bf w}(\bs)$.  The $\bu_{j_k}$ have prefixes that are either empty or of the form $\gws{\mm{s}{t'}}=\gws{i_s}$ ($s\not=t$). The primary difference conditions thereby prevent $i_t$-squares on $\gws{\mm{s}{t'}}\bu_{j_k}$.
The $\bu_{k}$ have as suffixes letters of the form $\gws{\mm{s}{t'}}$ for $1 < s \leq t' \leq r-1$ and $\gws{2i_1}$. The $i_t$-squarefreeness of $\bu_{j}$ guarantees the $i_t$-squarefreeness of $\bu_{j}\gws{2i_1}$. We must set conditions so that $\gw{1}{q-1}\gws{\mm{s}{t'}}\gw{1}{\mm{1}{t}}$ is also $\bs$-squarefree. In the case that $\mm{s}{t'} - i_t \leq i_t$ then the overlapping dissimilar pairs $i_t < \mm{s}{t'} \leq 2i_t < \mm{s}{t'} + i_t$  assure that $i_t$-squarefreeness continues over the concatenation $\bu_{k}\gws{\mm{s}{t'}}\bu_{{k+1}}$. This comes from $\gws{\mm{s}{t'}} = \gws{i_s}$ or $\gws{2i_t}$, $\gws{i_t} = \gws{\mm{s}{t-1} + i_t}$ and the primary difference conditions. 

Now consider the case $\mm{s}{t'} - i_r> i_r$. It is possible then that the primary difference conditions do not guarantee that the $i_t$-squarefreeness over $\bu_{k}\gws{\mm{s}{t'}}\bu_{{k+1}}$. (In particular, this can happen when $q_{k} \leq i_t$ together with some other conditions). We then set the secondary difference condition
\begin{align}\label{mu} \gws{i_s}=\gws{\mm{s}{t}} \not= \gws{\mm{s}{t}-i_k}.\end{align}
These conditions produce additional edges in the difference graph only if $\gws{\mm{s}{t}-i_k} \not=  \gws{2i_s}$ and $\gws{\mm{s}{t}-i_k} \not= \gws{i_j}$ for $j\not= s$. The colorability of the respective subgraph is maximized when the $\gws{\mm{s}{t}-i_k} = \gws{2i_f}$ for some $1<f\leq r$. Since the set of candidates $\gws{\mm{s}{t}}$ is at most $(r-2)$ candidates, namely $i_s$ for $1 < s < r$, these secondary difference conditions together with respective primary difference conditions produce a complete subgraph with the $r-1$ vertices $\{i_2, \ldots, i_{r-1}, 2i_f\}$ and this is $(r-1)$-colorable. Since this subgraph, together with the remainder of the difference graph, does not create a complete graph on $r+1$ vertices, the difference graph remains $r$-colorable.

{\bf II - } We now consider the general case. If there is a $2i_t \leq \mm{1}{t-1}$ for some $t$, then the action of $\tau_{i_t}$ replaces every subword $\gw{1}{\mm{1}{t-1}}$ in ${\bf w}(i_1, \ldots, i_{t-1})$ by
$$\gw{1}{i_t-1}\gws{i_t}\gw{1}{i_t-1}\gw{2i_t}{\mm{1}{t-1}}.$$ The action of all subsequent $\tau_{i_{t'}}$ ($t' > t$) will be restricted to the various copies of $\gw{1}{i_t}\gw{1}{i_t-1}\gws{2i_t}$ in $\bw(i_1,\ldots,i_t)$. The previous argument and conditions for squarefreeness applies to these subwords. 

Generally, the action of $\tau_{i_t}$ maps subwords of the form $$u\gws{\mm{s}{t}} = \gw{1}{\mm{1}{t-1}}^p\gw{1}{q-1}\gws{\mm{s}{t}}$$ (for some $p\geq 0, 0< q\leq \mm{1}{t-1},s,t$) onto:
$$\left[\gw{1}{i_t}\gw{1}{i-1}\gw{2i_t}{\mm{1}{t-1}}\right]^p\gw{1}{q-1}\gws{\mm{s}{t}}.$$
Squarefreeness requires $\gws{i_t} \not= \gws{2i_t}$.  
In addition, any prefixes of $\gw{1}{\mm{1}{t-1}}$ (that is, prefixes of $\bu$) will be mapped to prefixes of its image. By Corollary~\ref{cor2} and by induction, these prefixes are either the empty word or letters of type $\gws{\mm{s}{t}}= \gws{i_s}$ ($s \not= r$). So the primary difference conditions guarantee $\bs$-squarefreeness on those subwords of type $a\gw{1}{i_r-1}\gws{i_r}\gw{1}{i_r-1}\gws{2i_r}$ in $\bw(\bs)$ where $a = []$ or $\gws{\mm{s}{t}}$ (for some $s,t$).
 Conditions for subwords of type $\gw{1}{q-1}\gws{\mm{s}{t}}$ (and their concatenations with other subwords) were given in part I of this proof. We still need to impose secondary difference conditions on the subwords $\gw{2i_t}{\mm{1}{t-1}}$ because the primary difference conditions do not guarantee $i_t$-squarefreeness on them. There are possibly many ways in which this can be done in order to achieve this. We will give one way and then argue that the corresponding difference graph is $r$-colorable.

For
$k\geq t$ we denote by $\delta(t,k)$ the following set of difference conditions on the subwords $\bv:= \gw{2i_t}{\mm{1}{t-1}-1}$. Let $f_{t,k} := \max\{ f \ | \ 2i_t + fi_k < \mm{1}{t-1} \}$.
 Then in order that $\bv$ be $i_k$-squarefree, we set
\begin{align}\label{delta1} \gws{2i_t} \not= \gws{2i_t +i_k} \not= \cdots \not= \gws{2i_t+f_{t,k}i_k}.\end{align}
By induction, $\bv$ can occur as a prefix of $\bw_1:=\gws{\mm{1}{t-1}}\gw{1}{i_t}$ or, in the case that $\mm{1}{i_{t-1}}=i_{t-1}$ and $2i_{t-1} < \mm{1}{t-2}$, of $\bw_2:=\gw{2i_{t-1}}{\mm{1}{t-2}-1}$.  Although each $\bv$, $\bw_1$, and $\bw_2$ are $i_k$-squarefree, we need additional conditions in order to assure that their concatenations $\bv\bw_1$ and $\bv\bw_2$ are also $i_k$-squarefree. In particular, in $\bv\bw_1$, the conditions set up to this point provide  no guarantee that subwords of length $2i_k$ and containing $\gws{\mm{1}{t-1}}$ in the concatentation $\bv\bw_1$ are also $i_k$-squarefree. Therefore, we set
\begin{align}\label{delta2} \gws{\mm{1}{t-1}-1} \not= \begin{cases} \gws{\mm{1}{t-1}- i_k -1} & \text{ if } \mm{1}{t-1} - 2i_t > i_k \\
\gws{\mm{1}{t-1} - i_t - i_k -1} & \text{ otherwise, }
\end{cases}\end{align} the second case of (\ref{delta2}) taking into consideration when the word $\gw{2i_t}{\mm{1}{t-1}}$ is rather short.
The first thing to notice is that $\delta(t,k)$ is preserved by $\tau_{i_l}$ for all $l \geq t$ because $\tau_{i_l}$ fixes $\gw{2i_t}{\mm{1}{t-1}}$. We claim that the corresponding difference graph needs at most $r$ colors. In particular, we consider the ``worst case scenarios'' where the number of colors needed is maximized. When $\mm{1}{t-1} - 2i_t > i_t$ the most colors are needed when $t=2$, and every subsequent $i_{j}$ divides $i_{j-1}$. Condition (\ref{delta1}) of $\delta(2,2)$ produces a $2$-colorable subgraph. Adding the condition (\ref{delta1}) of $\delta(2,j)$ for every $j>2$ creates the need for one for additional color. The condition (\ref{delta2}) requires the most colors for the difference graph if every $\gws{\mm{1}{t-1}-i_k-1}$ is already a (distinct) node in the subgraph created by condition (\ref{delta1}). Then $\gws{\mm{1}{t-1}-1}$ would need to be a different color than these $r-1$ other nodes. The subsequent subgraph created by conditions (\ref{delta1}-\ref{delta2}) is therefore $r$-colorable. When $\mm{1}{t-1} - 2i_t \leq i_r$, the difference conditions come from the second case of (\ref{delta2}). The number of colors needed to color the difference graph is maximized when $\gws{\mm{1}{t-1}-1} = \gws{2i_2}$ and the $\gws{\mm{1}{t-1}-i_k-1}$ are all distinct and lie in $\{i_2,\ldots, i_r\}$. Together with the primary difference equations this creates a subgraph which is a complete subgraph on $r$ vertices, hence $r$-colorable. Again, since the ensuing difference graph coming from both primary and secondary difference conditions does not contain a complete graph on $r+1$ vertices. It is therefore $r$-colorable. All other cases produce difference subgraphs which are (possibly disjoint) unions of subgraphs of the above cases and therefore require fewer colors. \qed

\begin{Rem}\label{rem-ir1} When $i_r=1$, we can not use the same methods as above to determine $\bcr$ because $\tau_1$ acts trivially on $[2i_1]$. From Corollary~\ref{cor-2} and from Theorems~\ref{thm-equ},~\ref{thm-d} (to follow) we have that for $\bs = (1, i_{r-1}, \ldots, i_1)$, $\bcr(\bs) < \infty$ if and only if $\bs$ satisfies condition C. The difference graph is obtained by adding the following edges to the difference graph of $\bs' := (i_{r-1}, \ldots, i_1)$: 
$$\{\gws{s}, \gws{t} \} \mbox{ is an edge } \Leftrightarrow \exists x \in \forb{s}{i_{r-1}} \mbox{ such that } x+1 \in\forb{t}{i_{r-1}}.$$ For example, if $\bs' = (i_{r-1},2i_{r-1},\ldots, 2^{r-2}i_{r-1} )$, the corresponding difference graph is a complete graph on the $r-1$ vertices $i_{r-1},2i_{r-1},\ldots, 2^{r-2}i_{r-1}$. To obtain the difference graph for $\bs$, we add edges 
$ \{t,t+1\}$ for $1\leq t <i_{r-1}$, $\{1,i_t\}$ for $1 \leq t \leq r-1$,  $\{i_{r-1}-1, i_t\}$ for $1 \leq t \leq r-1$, and $\{i_{r-1}-1, 2i_1\}$. The resulting graph is $r+1$-colorable so $\bcr(1,i_{r-1},2i_{r-1},\ldots, 2^{r-2}i_{r-1} )=r+1$.

It is not hard to see that, generally, the vertex $\gws{1}$ in the difference graph has degree $r$ because it is connected to $\gws{2}$ and every $\gws{i_t}$, $i_t \in \bs$. All other vertices also have degree at most $r$. This gives $\bcr \leq r+1$. Unlike the case $i_r \geq 2$, the upper limit $r+1$ can be reached for sequences of type other than $(1, i_{r-1}, \ldots, 2^{r-2}i_{r-1})$: for example, $\bcr(1,2,5) = 4$ and $\bcr(1,3,5)=3$.
\end{Rem}

We now concentrate on sufficient conditions for $\bcr(\bs) = \infty$ to occur.

\begin{Satz}\label{thm-equ} Let $\bs= (i_r, \ldots, i_2,i_1) \in \BN_{\geq 1}^r$ ($r \geq 3$) be a strictly increasing sequence such that $i_1 = i_2 + \cdots i_{r-1} + i_r$. Then ${\bf cr}(\bs) = \infty$. \end{Satz}
{\bf Proof: } Without loss of generality, we can assume that $(i_1, \ldots,i_{r-1})$ and $(i_2, \ldots, i_r)$ satisfy condition C.  We first consider the case where $i_r \geq 2$. We show that $i_1 \sim i_r$. Since for $2 \leq j < r$, we have $i_j < i_j+ i_{j+1}+\cdots + i_r < 2i_j$, so $\tau_{i_j}(i_j + i_{j+1} + \cdots + i_r) = i_{j+1} + \cdots + i_r$. It follows:
\begin{eqnarray*}\tau_{i_1}\tau_{i_{r-1}}\tau_{i_{r-2}}\cdots \tau_{i_3}\tau_{i_2} (i_1) &=& \tau_{i_1}\tau_{i_{r-1}}\tau_{i_{r-2}}\cdots \tau_{i_3}\tau_{i_2} (i_2 + \cdots i_{r-1} + i_r) \\
&=& \tau_{i_1}\tau_{i_{r-1}}\tau_{i_{r-2}}\cdots \tau_{i_3}(i_3 + \cdots i_{r-1} + i_r)  = \cdots \\
\cdots &=& \tau_{i_1}(i_r) = i_r + i_1
\end{eqnarray*}
that is, $i_1 \sim i_r+i_1$. Since $i_r + i_1 \sim i_r$, we have $i_r \sim i_1$ and so ${\bf o}(\bs)$ has an $i_r$-square starting at $i_1$. 

Now suppose that $i_r = 1$. It suffices to show that $x \sim x+1$ for some $0 < x < 2i_1$ in ${\bf o}(i_1,\ldots, i_{r-1})$. But this follows from Corollary~\ref{cor-2}. \qed

For any $1 \leq s \leq t \leq r$ define $$x \mpl{\mm{s}{t}} := \begin{cases} x\mod \mm{s}{t} & \text{ if } x < \mm{s+1}{t} \text{ and } x\not\equiv 0\mod \mm{s}{t} \\ \mm{s}{t}  & \text{ if } x < \mm{s+1}{t} \text{ and } x\equiv 0 \mod \mm{s}{t} \\
x & \text{ otherwise } \end{cases}. $$

We say that a strictly increasing sequence $\bs = (i_r,\ldots, i_1) \in \BN^r_{\geq 2}$ satisfies {\bf condition D} if
\begin{eqnarray*} i_1 &<& i_2 + \cdots + i_{r-1}+  i_r  \\
 i_1 &>& i_2 + \cdots + i_{r-1} \\
i_2 &>&  i_3 + \cdots  + i_{r} \\
\vdots &> & \vdots \\
i_{r-2} &> & i_{r-1}  + i_r \\
i_{r-1} &> &i_r
\end{eqnarray*}
In other words, $(i_{r-1}, \ldots, i_{1})$ and $(i_r, \ldots, i_2)$ both satisfy condition C, but $(i_r, \ldots, i_1)$ doesn't.

\bigskip
\begin{Satz}\label{thm-d}  Let $\bs = (i_r, \ldots, i_2,i_1) \in \BN_{\geq 1}^r$ satisfy condition D.
Then $\bcr(\bs)=\infty$.
\end{Satz}
\noindent
{\bf Proof:} 
{\bf Case $i_r=1$:} The assumption $i_1-i_2 - \cdots - i_{r-1} < i_r = 1$ implies $i_1 = i_2 + \cdots + i_{r-1}$. So $\bcr(\bs)= \infty$ by Theorem~\ref{thm-equ}. 

{\bf Case $i_r \geq 2$:}
From Remark~\ref{rmk-1}, in order to prove that $\bcr(\bs)=\infty$ it suffises to show either that $i_r \sim 2i_r$ or that $i_r \sim i_j$ for some $j< r$ with $2i_r < i_j$. 

  Because $i_j < i_1 - i_2 - \cdots - i_{j-1} < 2i_j$ for every $1 < j < r$  we have $\mm{1}{j} = i_1 - i_2 - \cdots - i_j$ for every $1 < j <r $. From the definition of condition D we have that
$$\mm{1}{r-1} < i_r < \mm{2}{r-1} \leq \mm{3}{r-1} \leq \cdots \leq \mm{r-1}{r-1} < \mm{1}{r-2} < \mm{2}{r-2}$$
where $i_r < 2i_r <\mm{2}{r-2}$. The main difference between this case and that in Theorem~\ref{thm-main} is that here $\tau_{i_r}$ does {\it not} fix the orbit containing $i_r$. This is because $i_r \in \forb{i_r'}{i_{r-1}}$ where $i_r' := i_r \mpl{\mm{1}{r-1}} \leq \mm{1}{r-1} < i_r$. For every $0<x<i_r$, the action of $\tau_{i_r}$ on $x$ translates to the following action on ${\bf o}(i_1,\ldots, i_{r-1})$:
$$\forb{x \mpl{\mm{1}{r-1}}}{i_{r-1}} \act{i_j} \forb{f(x)}{i_{r-1}}$$
where $f(x) :=$ $$(\cdots (((x+i_r) \mpl{\mm{1}{r-2}} )\mpl{\mm{r-1}{r-1}}) \cdots \mpl{\mm{2}{r-1}}) \mpl{\mm{1}{r-1}}.$$
In particular, $\forb{x}{i_r}\supseteq \cup_k \forb{f^k(x)}{i_{r-1}}$. 

We now look at the action of $\tau_{i_r}$ on $i_r'$. If there is a $k > 0$ such that $f^k(i_r') = \mm{1}{r-2}$ or $= \mm{s}{r-1}$ for $s<r-1$ then we are done, because then $i_r \sim i_r' \sim i_s$ for some $s$ such that $2i_r < i_s$. If $f^k(i_r') < \mm{2}{r-1}$ for all $k$, then $$\{f^k(i_r') \ | \ k\geq 0\}\supseteq \{c, c^2, \cdots, \mm{1}{r-1}\}$$ where $c = \gcd(i_r, \mm{1}{r-1})$. We then have $\forb{i_r'}{i_r} \supset \forb{i_r'}{i_{r-1}} \cup \forb{\mm{1}{r-1}}{i_{r-1}}$ which gives $i_r \sim i_1$. In the case that $f^k(i_r') = \mm{r-1}{r-1} = i_{r-1}$ for some $k$, we have $2i_r>i_{r-1}$. So $$i_r \sim i_{r-1} \sim i_{r-1}-i_r\in \forb{(i_{r-1}-i_r)\mpl{\mm{1}{r-1}}}{i_{r-1}} = \forb{ld}{i_{r-1}}$$ where $d:=\gcd(i_{r-1}-i_r, \mm{1}{r-1}$ and $l \in \BN_{\geq 1}$. Without loss of generality, we can assume that $f^j(d)<\mm{2}{r-1}$ for all $j>0$ (otherwise we would be in the first case considered), so $\{f^k(i_r') \ | \ k\geq 0\}\supseteq \{d, d^2, \cdots, \mm{1}{r-1}\}$. It follows that $i_r \sim \mm{1}{r-1} \sim i_1$.

 The last case to consider is when $r=3$ and $2i_3 > i_1$. Although by the above argument $i_3 \sim i_1$, this does not give $i_3$-squareness starting at $i_1$ because $i_1 + 2i_3 > 2i_1$ and is therefore not defined. In this case we have $2i_3 \sim (2i_3-i_2)\mpl{(i_1-i_2)} = (i_3 \mpl{(i_2-i_3)})\mpl{(i_1-i_2)}$ so $2i_3 \sim i_3$ and we are done. 
\qed

\begin{Cor} Let $\bs \in \BN_{\geq 1}^r$ be strictly increasing and fail condition C. Then $\bcr(\bs) = \infty$. \end{Cor}
{\bf Proof:} If $\bs$ fails condition C, then it contains a subsequence $\bs' = (i_{j_k}, \ldots, i_{j_1})$ which satisfies either $i_{j_1} = i_{j_2} + \cdots + i_{j_k}$ or condition D. By Theorems~\ref{thm-equ} and~\ref{thm-d}, ${\bf o}(\bs) \preceq {\bf o}(\bs')$ must contain a square. \qed

Using all our previous results, we can now prove our main theorem:

{\bf Proof of Theorem~\ref{thm-main}:} Claim 1 is simply Lemma~\ref{lem-first}. Claim 2 is the result of Theorems~\ref{thm-pita},~\ref{thm-equ},~\ref{thm-d} and Remark~\ref{rem-ir1}. We note that $\bs$ fails Condition C if and only if there is a $t$ such that $i_t \leq i_{t+1} + \cdots + i_r$. To prove Claim 3 we need to consider possible dead-ends coming not only from $\bs$ but also any of its subsequences. A subsequence $\bs' < \bs$ creates no dead-end if $|\bs'| = l$ but $\bcr({\bf o}(\bs')) > l$. Since $\bcr({\bf o}(\bs'))$ either $= \infty$ or $ \leq |\bs'|+1$, we consider those subsequences of length $l$ and $\bcr \geq l+1$.  From Theorem~\ref{thm-pita} and Remark~\ref{rem-ir1}, this happens when $\bs'$ has the form $(1,k,2k,\ldots,2^{l-2}k)$, $(k,2k,\ldots,2^{l-1}k)$ or fails condition C. Then all longer subsequences $\bs' \subset \bs'' \subseteq \bs$ also satisfy $\bcr(\bs'') > l$. Therefore $G(\bs)$ has no dead-ends.
\qed

\begin{Rem} We conjecture that the conditions stated in Theorem~\ref{thm-main}, case {\bf 3)} are also necessary for $G(\bs)$ to have no dead-ends. Proving this, however, would require additional tools. The main problem is that when $r<l$, then any candidate for a dead-end ${\bf w}(i_{j_1},\ldots, i_{j_l})$ must also be $i_t$-square-free for all $i_t \in \bs$, $i_t \not\in (i_{j_1},\ldots, i_{j_l})$. A simple example is the following: let $A = \{ a, b\}$. Then $w:={\bf w}(4,5) = a^3ba^4ba$. If $\bs = (4,5)$, then $w$ would be a dead-end. However, if $\bs = (2,4,5)$, it cannot be because $w$ is not $2$-squarefree and would thereby not be a vertex of $G$. Despite this, $G$ does have dead-ends in this case. These come from ${\bf w}(2,5)=\gws{1}\gws{2}\gws{1}\gws{4}\gws{5}\gws{1}\gws{2}\gws{1}\gws{4}\gws{10}$. Setting $\gws{1}=\gws{4}=\gws{5}=a$ and $\gws{2}=\gws{10}=b$, we obtain $v = b^3a^2b^3ab$. Since $v$ is $2$, $4$, and $5$-squarefree, it is a vertex and, therefore, a dead-end of $G(2,4,5)$. 
\end{Rem}
\section{Conclusion}
Clearly these results are just the tip of the iceberg for this problem. In particular,  although Theorem~\ref{thm-main} tells us when we need to eliminate edges from the vertex set of $G(\bs)$ in order to obtain only non-trivial strongly connected components, it does tell us which ones to eliminate. Clearly, one needs to leave out the dead-ends. However, there are plenty of cases when a word is not a dead-end, but is incident only to dead-ends or squares. These must also then be eliminated. Last, considering the modelling problem on which this work is based, it would be interesting to extend this result to include other symmetries to which humans are sensitive, in particular mirror symmetries (palindromes) and, in the binary case, conjugation (word of type $w\tau(w)$ where $\tau$ permutes the two letters). 

\bibliography{paper}{}

\begin{thebibliography}{BEM79}

\bibitem[BEM79]{bean1979avoidable}
Dwight~R Bean, Andrzej Ehrenfeucht, and George~F McNulty.
\newblock Avoidable patterns in strings of symbols.
\newblock {\em Pacific J. Math}, 85(2):261--294, 1979.

\bibitem[Ber80]{Berstel1980235}
Jean Berstel.
\newblock Mots sans carr\'e et morphismes it\'er\'es.
\newblock {\em Discrete Mathematics}, 29(3):235 -- 244, 1980.

\bibitem[Ber84]{berstel1984some}
Jean Berstel.
\newblock Some recent results on squarefree words.
\newblock In {\em STACS 84}, pages 14--25. Springer, 1984.

\bibitem[Ber05]{berstel2005growth}
Jean Berstel.
\newblock Growth of repetition-free words \mbox{--} a review.
\newblock {\em Theoretical Computer Science}, 340(2):280--290, 2005.

\bibitem[BP07]{Berstel2007996}
Jean Berstel and Dominique Perrin.
\newblock The origins of combinatorics on words.
\newblock {\em European Journal of Combinatorics}, 28(3):996 -- 1022, 2007.

\bibitem[Bra88]{brandenburg1988uniformly}
Franz-Josef Brandenburg.
\newblock Uniformly growing ${\it k}^{th}$ power-free homomorphisms.
\newblock {\em Theoretical Computer Science}, 23(1):69--82, 1988.

\bibitem[Car83]{Carpi1983231}
Arturo Carpi.
\newblock On the size of a square-free morphism on a three letter alphabet.
\newblock {\em Information Processing Letters}, 16(5):231 -- 235, 1983.

\bibitem[Cro83]{Crochemore1983235}
M.~Crochemore.
\newblock Mots et morphismes sans carr\'e.
\newblock In C.~Berg\'e et~al., editor, {\em Combinatorial Mathematics
  Proceedings of the International Colloquium on Graph Theory and
  Combinatorics}, volume~75 of {\em North-Holland Mathematics Studies}, pages
  235 -- 245. North-Holland, 1983.

\bibitem[Cur93]{currie1993open}
James Currie.
\newblock Open problems in pattern avoidance.
\newblock {\em American Mathematical Monthly}, pages 790--793, 1993.

\bibitem[Cur05]{Currie20057}
James~D. Currie.
\newblock Pattern avoidance: themes and variations.
\newblock {\em Theoretical Computer Science}, 339(1):7 -- 18, 2005.
\newblock Combinatorics on Words.

\bibitem[Dej72]{Dejean197290}
Fran\c{c}oise Dejean.
\newblock Sur un th\'eor\`eme de {T}hue.
\newblock {\em Journal of Combinatorial Theory, Series A}, 13(1):90 -- 99,
  1972.

\bibitem[FK97]{FK1}
Ruma Falk and Clifford Konold.
\newblock Making sense of randomness: Implicit encoding as a basis for
  judgment.
\newblock {\em Psychol Rev}, 104(2):301--318, 1997.

\bibitem[GT04]{GT3}
Thomas Griffith and Joshua Tenenbaum.
\newblock From algorithmic to subjective randomness.
\newblock In {\em Advances in Neural Information Processing Systems},
  volume~16, pages 953--960, 2004.

\bibitem[Lot97]{lothaire1997combinatorics}
M\_ Lothaire.
\newblock {\em Combinatorics on words}.
\newblock Cambridge University Press, 1997.

\bibitem[LVH06]{Leong:2006:RRD:1142405.1142428}
Tuck~Wah Leong, Frank Vetere, and Steve Howard.
\newblock Randomness as a resource for design.
\newblock In {\em Proceedings of the 6th conference on Designing Interactive
  systems}, DIS '06, pages 132--139, New York, NY, USA, 2006. ACM.

\bibitem[NZ99]{noonan1999goulden}
John Noonan and Doron Zeilberger.
\newblock The {G}oulden \mbox{—-} {J}ackson cluster method: extensions,
  applications and implementations.
\newblock {\em Journal of Difference Equations and Applications},
  5(4-5):355--377, 1999.

\bibitem[Sha04]{shallit2004simultaneous}
Jeffrey Shallit.
\newblock Simultaneous avoidance of large squares and fractional powers in
  infinite binary words.
\newblock {\em International Journal of Foundations of Computer Science},
  15(02):317--327, 2004.

\bibitem[Thu12]{thue1912gegenseitige}
Axel Thue.
\newblock {\em {\"U}ber die gegenseitige Lage gleicher Teile gewisser
  Zeichenreihen}.
\newblock J. Dybwad, 1912.

\end{thebibliography}
\bibliographystyle{alpha}

\bigskip
\leftline{Yasmine B. Sanderson}

\leftline{Emmy-Noether-Zentrum}

\leftline{Friedrich-Alexander-Universit\"{a}t}

\leftline{Cauerstra\ss e 11}

\leftline{91054 Erlangen, Germany}
 
\leftline{Email: sanderson@mi.uni-erlangen.edu}
\end{document}